\documentstyle[12pt]{article}

\begin{document}
\newtheorem{izrek}{Theorem}[section]
\newtheorem{posledica}[izrek]{Corollary}
\newtheorem{lema}[izrek]{Lemma}
\newtheorem{trditev}[izrek]{Proposition}

\newcommand{\MM}{M^{\ast}}
\newcommand{\En}{{\cal E}\!\ell_0}
\newcommand{\REn}{{\cal RE}\!\ell_0}
\newcommand{\JEn}{{\cal J}\!{\cal E}\!\ell_0}
\newcommand{\E}{{\cal E}\!\ell}
\newcommand{\RE}{{\cal RE}\!\ell}
\newcommand{\JE}{{\cal J}\!{\cal E}\!\ell}
\newcommand{\A}{{\cal A}}
\newcommand{\Ap}{{\cal A}^{\prime}}
\newcommand{\B}{{\cal B}}
\newcommand{\F}{F[X]}
\newcommand{\D}{{\cal D}}
\newcommand{\DD}{{\cal D}^{\prime}}

\title{LOCAL AUTOMORPHISMS OF THE UNITARY GROUP AND THE GENERAL LINEAR
GROUP ON A HILBERT SPACE}
\author{Lajos Moln\' ar\\
        Institute of Mathematics\\
        Lajos Kossuth University\\
        P.O. Box 12\\
        4010 Debrecen, Hungary\\
        \texttt{e-mail:molnarl@math.klte.hu}
        \\[5mm]
        Peter \v Semrl\\
        Department of Mathematics\\
        University of Ljubljana\\
        Jadranska 19\\
        1000 Ljubljana, Slovenia\\
        \texttt{e-mail:peter.semrl@fmf.uni-lj.si}}
\date{}
\maketitle
\newpage
\vglue 8truecm
\begin{abstract}
\noindent
We prove that every 2-local automorphism of the unitary group or the
general linear group on a complex infinite-dimensional separable Hilbert
space is an automorphism. Thus these types of transformations
are completely determined by their local actions on the two-points
subsets of the groups in question. \end{abstract}

\section{Introduction}
The study of automorphism groups of algebraic structures is of
great importance in every field of mathematics.
In a series of papers (see \cite{BaMo, MoGyo, Mol2} and the references
therein) we investigated these groups from the
point of view of how they are determined by their local actions.
Our investigations were motivated by the paper \cite{Kad} of Kadison on
local derivations and by a problem of Larson in \cite{Lar}
initiating the study of local automorphisms of Banach
algebras. The structures
that we treated so far were mainly $C^*$-algebras and we considered the
following question: When is it true that any local automorphism, that
is,
any linear transformation which pointwise equals an automorphism (this
automorphism may, of course, differ from point to point) is an
automorphism?

It is easy to see that if we drop the assumption of the linearity of the
transformations
in question, then the corresponding statements are no longer true.
However, if instead of linearity plus locality we assume the so-called
2-locality, then we can obtain positive results (for the first such
result see \cite{Se2}).
2-locality means that our transformation (linearity is not assumed any
more) is supposed to be equal to an automorphism at every pair of
points.
Notice that in this way we arrive at a question that can be raised in
any algebraic structure. For example, this observation motivated us to
consider the problem for
the orthomodular poset of all projections on a Hilbert space whose
structure plays a fundamental role in the mathematical foundations of
quantum mechanics (see \cite{Mol1}).

In the present paper we study the analogous problem in the case of two
important groups appearing in pure algebra and in the theory of
operator algebras. They are the unitary group and the general linear
group.

We begin with some notation and definitions.
Throughout the paper $H$ denotes a complex infinite-dimensional
separable
Hilbert space. By $B(H)$, $U(H)$, $GL(H)$, and $A(H)$ we denote the
algebra of all bounded linear operators on $H$, the multiplicative group
of all unitary operators on $H$, the general linear group on $H$ consisting
of all bounded invertible linear operators on $H$, and the set of all
additive maps on $H$ (that is, the set of all maps $A:H\to H$ satisfying
$A(x+y)=Ax + Ay$ $(x,y\in H)$), respectively. For $T\in B(H)$ we denote
its spectrum by $\sigma (T)$.

A mapping $\phi : GL(H) \to GL(H)$ is called a 2-local automorphism of
the general linear group if for every $X,Y\in GL(H)$ there is an
automorphism $\phi_{X,Y}$ of the group $GL(H)$, depending on
$X$ and $Y$,
such that $\phi_{X,Y} (X) = \phi (X)$ and $\phi_{X,Y} (Y) = \phi (Y)$.
In the case of the unitary group the situation is somewhat different.
Because of certain reasons (see, for example, Theorem~\ref{I:Sakai}),
$U(H)$ is considered
here as a topological group and by an automorphism of $U(H)$ we mean a
uniformly continuous group-automorphism.
Now, not surprisingly,
a mapping $\phi : U(H) \to U(H)$ is called a 2-local automorphism
of the unitary group if for every $X,Y\in U(H)$ there is an
automorphism  $\phi_{X,Y}$ of $U(H)$ (in the above sense),
such that $\phi_{X,Y} (X) = \phi (X)$ and $\phi_{X,Y} (Y) = \phi (Y)$.

\section{Local automorphisms of the unitary group}

It was proved in \cite{Sak} that any uniformly continuous group
isomorphism between the unitary groups of two $AW^*$-factors is
implemented by a linear or conjugate-linear *-isomorphism of the
factors themselves. As a particular case, concerning $B(H)$ we have the
following result.

\begin{izrek}\label{I:Sakai}
{\rm {(Sakai)}} Let $\phi: U(H) \to U(H)$ be a uniformly continuous
automorphism. Then there exists either a unitary or an antiunitary
operator $U$ on $H$ such that
\begin{equation}\label{E:elso}
\phi(V)=UVU^* \qquad (V\in U(H)).
\end{equation}
\end{izrek}

As for the 2-local automorphisms of the unitary group $U(H)$ we have the
following statement. In the proof we use the notation $x\otimes y$ which
denotes the operator on $H$ defined by
$$
(x\otimes y)(z)=\langle z,y\rangle x \qquad (z\in H)
$$
for any $x,y\in H$.

\begin{izrek}\label{uH}
Every 2-local automorphism of $U(H)$ is an automorphism.
\end{izrek}

{\bf Proof}.
Let $\phi:U(H) \to U(H)$ be a 2-local automorphism.
For every projection $P\in B(H)$, the operator $I-2P$ is unitary.
Since $\phi$ is locally of the form $(\ref{E:elso})$, we obtain that
$\phi(I-2P)=I-2P'$ for some projection $P'$. Consider the transformation
$$
P \longmapsto (I-\phi(I-2P))/2.
$$
This is a 2-local automorphism of the orthomodular poset of all
projections on $H$. By \cite[Proposition]{Mol1}, it is an automorphism
and there is either a unitary or an antiunitary operator $U$ on $H$
such that our transformation is of the form
$$
P \longmapsto (I-\phi(I-2P))/2=UPU^*.
$$
We have $\phi(I-2P)=U(I-2P)U^*$ for every projection $P$.
Transforming the original map $\phi$ by this operator $U$, we can
suppose without loss of generality that
$\phi(I-2P)=I-2P$ for every projection $P$.
Let $V\in U(H)$ and pick an
arbitrary unit vector $x\in H$.
Let $P$ be the orthogonal projection onto the subspace spanned by $x$,
that is, let $P=x\otimes x$. By the local
property of $\phi$ we have an either unitary or antiunitary operator
$U_{V,P}$ such that
$\phi(V)=U_{V,P}VU_{V,P}^*$ and $\phi(I-2P)=U_{V,P}(I-2P)U_{V,P}^*$.
Since $\phi(I-2P)=I-2P$, it follows that $P=U_{V,P}PU_{V,P}^*$.
We compute
$$
\langle \phi(V) x,x \rangle x \otimes x=
x \otimes x \cdot \phi(V) \cdot x \otimes x=
U_{V,P} PU_{V,P}^* U_{V,P} VU_{V,P}^* U_{V,P} PU_{V,P}^*=
$$
$$
U_{V,P} P V PU_{V,P}^*=
U_{V,P} \cdot \langle Vx,x \rangle x\otimes  x \cdot U_{V,P}^*.
$$
Since $U_{V,P}$ is either linear or conjugate-linear, we have
either
$$
\langle \phi(V) x,x \rangle x \otimes x=
\langle Vx,x \rangle U_{V,P} x\otimes  U_{V,P} x
$$
or
$$
\langle \phi(V) x,x \rangle x \otimes x=
{\overline{\langle Vx,x \rangle}} U_{V,P} x\otimes  U_{V,P} x.
$$
We deduce that for every $x\in H$ either
$$
\langle \phi(V) x,x \rangle=
\langle Vx,x \rangle
$$
or
$$
\langle \phi(V) x,x \rangle=
{\overline{\langle Vx,x \rangle}}
$$
holds true.
It is rather elementary to verify (see \cite[Lemma]{Mol1}) that this
implies that either $\phi(V)=V$ or $\phi(V)=V^*$.
We show that either
$\phi(V)=V$ for every $V\in U(H)$ or $\phi(V)=V^*$ for every $V\in
U(H)$. To see this, observe that $\phi(iI)$ is either $iI$ or $-iI$.
First assume that $\phi(iI)=iI$. We assert that
in this case we have $\phi(V)=V$ $(V\in U(H))$. Suppose on the contrary
that
there is a non-selfadjoint unitary operator $V$ for which $\phi(V)=V^*$.
Let $\lambda \in \sigma(V)$. By the spectral theorem of normal operators
we can choose a sequence $(x_n)$ of pairwise orthogonal unit vectors
in $H$
such that $\langle Vx_n, x_n\rangle \to \lambda$. We extend $(x_n)$ to a
complete orthonormal sequence $(x_n')$ in $H$. Pick pairwise different
complex
numbers $\lambda_n$ of modulus 1 from the open upper half-plane and
consider
the unitary operator $W=\sum_n \lambda_n x_n'\otimes x_n'$. By the local
property of $\phi$ we have an either unitary or antiunitary operator
$U_{i,W}$ such that $\phi(iI)=U_{i,W} iI U_{i,W}^*$ and
$\phi(W)=U_{i,W}WU_{i,W}^*$.
Since we have supposed that $\phi(iI)=iI$, it follows that $U_{i,W}$ is
unitary.
So, $\phi(W)= \sum_n \lambda_n U_{i,W}x_n'\otimes U_{i,W}x_n'$ and,
on the other hand, we know that $\phi(W)=W$ or $\phi(W)=W^*$. These
result in $\phi(W)=W$.
Once again, by the local property of $\phi$ we have an either unitary or
antiunitary operator $U_{W,V}$ such that
$$
\phi(W)= U_{W,V} W U_{W,V}^* \quad {\rm{ and }} \quad
\phi(V)= U_{W,V} V U_{W,V}^*.
$$
Since $\phi(W)=W$, it follows that $U_{W,V}$ is necessarily linear and
from the equalities
$$
\sum_n \lambda_n x_n'\otimes x_n'=W=
\phi(W)= U_{W,V} W U_{W,V}^*=
\sum_n \lambda_n U_{W,V}x_n'\otimes U_{W,V}x_n'
$$
we conclude that $U_{W,V}$ is diagonizable with respect to $(x_n')$.
Therefore, we can compute
$$
\langle x_n , Vx_n\rangle =
\langle V^* x_n , x_n\rangle =
\langle \phi(V) x_n , x_n\rangle =
$$
$$
\langle U_{W,V} V U_{W,V}^* x_n , x_n\rangle =
\langle V U_{W,V}^* x_n , U_{W,V}^* x_n\rangle =
\langle V x_n , x_n\rangle .
$$
If $n$ goes to infinity, we obtain that $\lambda=\overline{\lambda}$.
Since $\lambda $ was an arbitrary element of $\sigma(V)$, we infer that
$V=V^*$ which is a contradiction. So, we have $\phi(V)=V$ for every
$V\in U(H)$ and this shows that $\phi$ is an automorphism of $U(H)$.
We now consider the case when $\phi(iI)=-iI$. Similarly as above one can
verify that then we have $\phi(V)=V^*$ $(V\in U(H))$. By the local
property of $\phi$ it follows that for every $V,V'\in U(H)$ there exists
an either unitary or antiunitary operator $U$ such that
$V^*=UVU^*$ and ${V'}^*=UV'U^*$. If $V$ is diagonal
with respect to an orthonormal basis defined in the same way as $W$
and $V'$ permutes the same basis, then one can easily arrive at a
contradiction.

The proof of the theorem is now complete.

\section{Local automorphisms of the general linear\\ group}

We first remark that we were unable to
find the description of the general form of automorphisms of $GL(H)$
in the literature. However, applying a result of Radjavi \cite{Rad}
on a factorization of invertible operators into a product of involutions
and some automatic continuity techniques it is possible to obtain the
general form of the automorphisms of $GL(H)$ as a consequence of
results
of Rickart on the isomorphisms of some analogues of the classical groups
\cite{Ri1,Ri2}. So, we begin with a statement on the form of the
automorphisms of the general linear group $GL(H)$.

\begin{izrek}\label{formglH}
Let $\phi: GL(H)\to GL(H)$ be an automorphism of the general linear
group. Then $\phi$ is of one of the following forms:
\begin{itemize}
\item[(i)]
there exists a bounded linear invertible operator $T:H\to H$ such that
$$\phi (X) = TXT^{-1} \qquad (X\in GL(H)),$$

\item[(ii)]
there exists a bounded conjugate-linear invertible operator $T:H\to H$ such that
$$\phi (X) = TXT^{-1} \qquad (X\in GL(H)),$$

\item[(iii)]
there exists a bounded linear invertible operator $T:H\to H$ such that
$$\phi (X) = \left(TX^{-1}T^{-1}\right)^* \qquad (X\in GL(H)),$$

\item[(iv)]
there exists a bounded conjugate-linear invertible operator $T:H\to H$ such that
$$\phi (X) = \left(TX^{-1}T^{-1}\right)^* \qquad (X\in GL(H)).$$
\end{itemize}
\end{izrek}
\noindent \bf Remark. \rm
If $\phi$ is of type (i), (ii), (iii), (iv), then for every
$X\in GL(H)$ we have $\sigma (\phi (X))=\sigma (X)$,
$\sigma (\phi (X))=\overline{\sigma (X)}$, $\sigma (\phi (X))=(\overline{\sigma (X)})^{-1}$,
$\sigma (\phi (X))=\sigma (X)^{-1}$, respectively.   \\

{\bf Proof}.
Let $\phi: GL(H)\to GL(H)$ be an automorphism.
By a result of Rickart \cite[Theorem 5.1]{Ri1}, \cite[Theorem I]{Ri2},
$\phi$ must be either of the form $\phi(X)=\tau (X) TXT^{-1}$ or of the form
$\phi(X)=\tau(X) (TX^{-1}T^{-1})^*$, where $\tau: GL(H) \to {\bf C}$ is a multiplicative
map and $T:H\to H$ is a bijective additive map satisfying $T(\lambda x) = f(\lambda )Tx$,
$\lambda\in {\bf C}$, $x\in H$, for some ring automorphism $f$ of {\bf C}
(one has to be careful when applying the result of Rickart since $S^*$ in his papers
denotes the adjoint of an operator $S$ defined as for Banach space operators, while here,
of course, $S^*$ denotes the adjoint in the Hilbert space sense).
Let us first show that
$\tau (X)\equiv 1$. Since $\tau (I) = 1$ we have $\tau (S) \in \{ -1,1 \}$ for any involution
$S\in GL(H)$, that is, for any $S$ satisfying $S^2 = I$. According to \cite{Rad} every
element of $GL(H)$ can be written as a product of involutions, and consequently, the range
of $\tau$ is contained in $\{ -1,1 \}$. An arbitrary involution $S$ can be expressed as
$S= (I-P)-P$ where $P$ is an idempotent. For $R= (I-P) + iP$ we have $S=R^2$, and therefore,
$\tau (S) = 1$. Applying \cite{Rad} once again we conclude that $\tau$ is identically
equal to 1.

Therefore, either $\phi(X)= TXT^{-1}$ or
$\phi(X)= (TX^{-1}T^{-1})^*$. In the second case we can compose $\phi$ by
$X\mapsto (X^* )^{-1}$ to conclude that in both cases $X\mapsto TXT^{-1}$ is
an automorphism of $GL(H)$. We have to prove that $f:{\bf C}\to {\bf C}$ is
either the identity or the complex conjugation and that $T$ is bounded.
To prove this one can apply automatic continuity techniques (note that for
this part of the proof the assumption that $H$ is infinite-dimensional is
indispensable). However, we shall use a shorter way of reducing the
problem to a known result.

So, assume that $\phi(X)= TXT^{-1}$ is an automorphism of $GL(H)$. Clearly,
$\phi^{-1}(X)= T^{-1}XT$. Define additive mappings
$\psi, \varphi : B(H)\to A(H)$ by
$\psi(X)= TXT^{-1}$ and
$\varphi(X)= T^{-1}XT$, $X\in B(H)$. If $|\lambda | > || X ||$, then
$\psi (X) = \psi (X - \lambda I) + \psi (\lambda I) \in GL(H) + GL(H)\subset
B(H)$. So, $\psi$ is a multiplicative map from $B(H)$ into $B(H)$ which
is also bijective since $\varphi: B(H) \to B(H)$ is its inverse. By
\cite{Se1} there exists a bijective bounded linear or conjugate-linear
map $S:H\to H$ such that $TXT^{-1} = SXS^{-1}$ for every $X\in B(H)$,
or equivalently, the additive map $S^{-1}T$ commutes with every $X\in B(H)$.
It follows that $T=\lambda S$ for some nonzero scalar $\lambda$.
This completes the proof of the statement that every automorphism
$\phi$ of $GL(H)$ has one of the forms (i), (ii), (iii) or (iv).\\

In the proof of the main result of this section we shall need the
following lemma.
Let $K$ be a nonempty subset of the complex plane and let $\lambda$
be a complex number. We use the following notation:
$K-\lambda = \{ \mu - \lambda \, : \, \mu \in K \}$,
$\overline{K} = \{ \overline{\mu} \, : \, \mu \in K \}$,
and $r(K) = \sup \{ | \mu | \, : \, \mu\in K \}\in [0, \infty ]$. If also $0\not\in K$,
then $K^{-1} = \{ \mu^{-1}  \, : \, \mu \in K \}$.\\

\noindent
{\bf Lemma.}
{\it Let $K$ be a nonempty
compact subset of {\bf C} and $\lambda$ a complex number such that
${\rm Im}\, \lambda > r(K)+1$. Then $0\not\in K - \lambda$, $K-\lambda \not=
(K - \lambda )^{-1}$, $K-\lambda \not= \overline{K - \lambda}$, and
$K-\lambda\not= (\overline{K-\lambda})^{-1}$.}\\

{\bf Proof}. Clearly, $0\not\in K - \lambda$ and $|\lambda | > r(K) + 1$. It follows
that $r(K-\lambda) \ge |\lambda | - r(K) > 1$. On the other hand,
$r((K - \lambda)^{-1}) = { 1 \over |\lambda - \lambda_0 |}$ for some $\lambda_0
\in K$. Since
$$
{1 \over |\lambda - \lambda_0 |} \le
{1 \over |\lambda | - |\lambda_0 |} \le
{1 \over |\lambda | - r(K)}
<1,
$$
we have $r(K-\lambda ) >  r((K-\lambda)^{-1}) = r((\overline{K-\lambda})^{-1})$,
which further yields that $K-\lambda \not= (K-\lambda)^{-1}$ and
$K-\lambda \not= (\overline{K-\lambda})^{-1}$. Furthermore, our assumption implies
that $K-\lambda$ belongs to the open lower half-plane, and consequently,
$K-\lambda\not= \overline{K-\lambda}$. This completes the proof.

\begin{izrek}\label{glH}
Every 2-local automorphism of $GL(H)$ is an automorphism.
\end{izrek}

{\bf Proof}.
Assume that $\phi$ is a 2-local automorphism of $GL(H)$. Composing
it with an appropriate automorphism of $GL(H)$ we can assume with no
loss of generality that $\phi (2iI) = 2iI$. It follows from the Remark
that $\phi_{X, 2iI}$ has to be of type (i) for every $X\in GL(H)$.
In particular, we have $\sigma (\phi (X))=\sigma (X)$ $(X\in GL(H))$,
and $\phi(\lambda I)= \lambda I$ $(\lambda \in {\bf C})$. Denote by
${\cal S}$ the set of all operators $X\in GL(H)$ satisfying
$\sigma (X)\not=\overline{\sigma (X)}$ and $\sigma (X)\not=\sigma (X)^{-1}$,
$\sigma (X)\not=(\overline{\sigma (X)})^{-1}$. If $X\in {\cal S}$ and $Y$
is an arbitrary element of $GL(H)$ then $\phi_{X,Y}$ has to be of type (i).
In particular, $\phi_{X,Y} (\lambda I) = \lambda I$ $(\lambda\in {\bf
C})$.

Now, let $X\in B(H)$ be any bounded linear operator on $H$. Denote by
$L_X$ the set of all complex numbers $\lambda$ such that $X-\lambda I$
is an invertible operator contained in ${\cal S}$. By Lemma, this
set is always nonempty. We define $\psi: B(H)\to B(H)$ by
$\psi (X) = \phi (X-\lambda I) + \lambda I$ where $\lambda\in L_X$.
First we have to show that $\psi$ is well-defined. So, assume that
$\mu$ also belongs to $L_X$. Then we already know that
$\phi_{ X-\lambda I, X -\mu I}$
is of type (i), and consequently,
\[
\phi (X-\lambda I) + \lambda I =
\phi_{ X-\lambda I, X -\mu I}(X-\lambda I) + \lambda I =
\]
\[
\phi_{ X-\lambda I, X -\mu I}(X) = \phi (X-\mu I) + \mu I.
\]
Next, we
show that the restriction of $\psi$ to $GL(H)$ coincides with $\phi$.
In order to do this we first observe that if $\varphi : GL(H) \to GL(H)$
is an automorphism of type (i) and if $X,Y$ are arbitrary elements of
$GL(H)$ such that $X+Y$ is also invertible, then $\varphi (X+Y) =
\varphi (X) + \varphi (Y)$. So, for any $X\in GL(H)$ and $\lambda
\in L_X \setminus \{ 0 \}$ we have
\[
\psi (X) = \phi (X-\lambda I) +\lambda I = \phi_{X, X-\lambda I}
(X-\lambda I) + \lambda I =
\]
\[
\phi_{X, X-\lambda I}(X) + \phi_{X, X-\lambda I}(-\lambda I)
+\lambda I = \phi_{X, X-\lambda I}(X) = \phi (X).
\]
It is well-known that every algebra automorphism of $B(H)$ is
inner. We show that $\psi$ is a 2-local automorphism of $B(H)$,
that is, for every pair $X,Y\in B(H)$ there exists a bounded linear invertible
$S:H\to H$ such that $\psi (X) = SXS^{-1}$ and
$\psi (Y) = SYS^{-1}$. By Lemma we know that there exists $\lambda\in
{\bf C}$ such that $\lambda\in L_X \cap L_Y$. The automorphism
$\phi_{X-\lambda I, Y-\lambda I}$
is of type (i), and therefore spatially implemented by a bounded linear
invertible operator, say $S$. Then
\[
\psi (X) = \phi (X-\lambda I) + \lambda I
= \phi_{X-\lambda I, Y-\lambda I} (X-\lambda I) + \lambda I =
\]
\[
S(X-\lambda I)S^{-1} + \lambda I = SXS^{-1},
\]
and similarly,
\[
\psi (Y) = SYS^{-1}.
\]

Applying the result of the second author \cite{Se2} on 2-local automorphisms of
$B(H)$ we conclude that there exists a bounded linear invertible operator $T:
H\to H$ such that $\psi (X) = TXT^{-1}$ for every $X\in B(H)$.
Consequently, $\phi (X) = TXT^{-1}$ for every $X\in GL(H)$.
This completes the proof.\\

\section*{Acknowledgement}

We express our gratitude to Professor R.V.
Kadison for providing information on the mentioned papers of
C.E. Rickart.

\noindent
This research was supported by a Hungarian-Slovene grant and also by a
grant from the Ministry of Science of Slovenia.

\baselineskip16pt
{\small

}

\vskip 1cm
\noindent
Lajos Moln\' ar\\
Institute of Mathematics\\
Lajos Kossuth University\\
P.O. Box 12\\
4010 Debrecen, Hungary\\
\texttt{e-mail:molnarl@math.klte.hu}        \\[5mm]
Peter \v Semrl\\
Department of Mathematics\\
University of Ljubljana\\
Jadranska 19\\
1000 Ljubljana, Slovenia\\
\texttt{e-mail:peter.semrl@fmf.uni-lj.si}


\begin{thebibliography}{99}

\bibitem{BaMo}
C.J.K. Batty and L. Moln\'ar,
{\em On topological reflexivity of the groups of
\linebreak *-automorphisms and surjective isometries of $\mathcal
B(H)$,} Arch. Math. \textbf{67} (1996), 415--421.

\bibitem{Kad}
R.V. Kadison,
{\em Local derivations},
J. Algebra \textbf{130} (1990), 494--509.

\bibitem{Lar}
D.R. Larson,
{\em Reflexivity, algebraic reflexivity and linear interpolation},
Amer. J. Math. \textbf{110} (1988), 283--299.

\bibitem{Mol2}
L. Moln\'ar,
{\em Reflexivity of the automorphism and isometry groups of
$C^*$-algebras in BDF theory,}
Arch. Math. \textbf{74} (2000), 120--128.

\bibitem{Mol1} L. Moln\'ar, {\em Local automorphisms of some quantum
mechanical structures}, preprint

\bibitem{MoGyo}
L. Moln\'ar and M. Gy\H{o}ry,
{\em Reflexivity of the automorphism and isometry groups of the
suspension of $\mathcal B(\mathcal H)$,}
J. Funct. Anal. \textbf{159} (1998), 568--586.


\bibitem{Rad} H. Radjavi, {\em The group generated by involutions},
Proc. Roy. Irish Acad. Sect. A {\bf 81} (1981), 9--12.

\bibitem{Ri1} C.E. Rickart, {\em Isomorphic groups of linear
transformations}, Amer. J. Math. {\bf 72} (1950), 451--464.

\bibitem{Ri2} C.E. Rickart, {\em Isomorphisms of infinite-dimensional
analogues of the classical groups}, Bull. Amer. Math. Soc. {\bf 57}
(1951), 435--448.

\bibitem{Sak} S. Sakai, {\em On the group isomorphisms of unitary groups
in $AW^*$-algebras}, T\^ ohoku Math. J. {\bf 7} (1955), 87--95.

\bibitem{Se1} P. \v Semrl, {\em Isomorphisms of standard operator algebras},
Proc. Amer. Math. Soc. {\bf 123} (1995), 1851--1855.

\bibitem{Se2} P. \v Semrl, {\em Local automorphisms and derivations on ${\cal B}(H)$},
Proc. Amer. Math. Soc. {\bf 125} (1997), 2677--2680.


\end{thebibliography}
\end{document}